\documentclass[10pt]{article}

\usepackage{amsmath}
\usepackage{amssymb}
\usepackage{eucal}

\begin{document}

\baselineskip 16pt

\title{ON A FINITE GROUP HAVING A NORMAL SERIES WHOSE
FACTORS HAVE BICYCLIC SYLOW SUBGROUPS}
\author{V.\,S.\,Monakhov and A.\,A.\,Trofimuk}

\date{}
\maketitle
\begin{abstract}

We consider the structure of a finite groups having a normal
series whose factors have bicyclic Sylow subgroups. In particular,
we investigated groups of odd order and $A_4$-free groups with
this property. Exact estimations of the derived length and
nilpotent length of such groups are obtained.

\end{abstract}

{\small       {\bf Keywords}: normal series, bicyclic Sylow
subgroups, $A_4$-free groups, derived length, nilpotent length.}

{\small }

\section{INTRODUCTION}

All groups considered in this article will be finite.

By the Zassenhaus Theorem (see Huppert, 1967, IV, 2.11) the
derived subgroup of a group with cyclic Sylow subgroups is a
cyclic Hall subgroup such that the corresponding quotient group is
also cyclic. Hence the derived length of such group is at most 2.

Recall that a group is bicyclic if it is the product of two cyclic
subgroups. The invariants of the groups with bicyclic Sylow
subgroups were found in Monakhov, Gribovskaya (2001). In
particular, it is proved that the derived length of such groups is
at most 6 and the nilpotent length of such groups is at most 4.

Let the group $G$ have a normal series in which every Sylow
subgroup of its factors is cyclic. Then $G$ is supersolvable by
the Zassenhaus Theorem.

In this article we study groups having a normal series whose
factors have bicyclic Sylow subgroups. We prove the following

{\bf Theorem 1.1.} {\sl Let $G$ be a solvable group having a
normal series such that every Sylow subgroup of its factors is
bicyclic. Then the following statements hold:

$1)$ the nilpotent length of $G$ is at most 4 and the derived
length of $G/\Phi (G)$ is at most~5;

$2)$ $G$ contains a normal subgroup $N$ such that $G/N$ is
supersolvable and $N$ possesses an ordered Sylow tower of
supersolvable type;

$3)$ $l_{2}(G) \le 2$, $l_{3}(G) \le 2$ and $l_{p}(G) \le 1$ for
every prime $p>3$;

$4)$ $G$ contains a normal Hall
$\{2,3,7\}^{\prime}$\nobreakdash-\hspace{0pt}subgroup $H$ and $H$
possesses an ordered Sylow tower of supersolvable type.}

Here $\Phi (G)$ is the Frattini subgroup of $G$ and $l_{p}(G)$ is
the $p$-length of $G$. A group $G$ is $A_4$-free if there is no
section isomorphic to the alternating group $A_4$ of degree 4.

{\bf Corollary 1.2.} {\sl  Let $G$ be a solvable group having a
normal series such that every Sylow subgroup of its factors is
bicyclic. If $G$ is an $A_4$-free group then the following
statements hold:

$1)$ $l_{p}(G) \le 1$ for every prime $p$;

$2)$ the derived length of $G/\Phi (G)$ is at most~3.}

{\bf Corollary 1.3.} {\sl  Let $G$ be a group of odd order having
a normal series such that every Sylow subgroup of its factors is
bicyclic. Then the following statements hold:

$1)$ $G$ possesses an ordered Sylow tower of supersolvable type;

$2)$ the derived subgroup of $G$ is nilpotent. In particular,
$G/\Phi (G)$ is metabelian.}

Examples that show accuracy of the estimations in Theorem 1.1 and
Corollary 1.2 are constructed, see Examples 3.1 -- 3.3.

\section{PRELIMINARIES}

In this section, we give some definitions and basic results which
are essential in the sequel.

A normal series of a group $G$ is a finite sequence of  normal
subgroups $G_i$ such that

$$
1=G_0\subseteq G_1\subseteq \ldots \subseteq G_m=G. \eqno (1)
$$

We call the groups $G_{i+1}/G_i$ the factors of the normal series
$(1)$.

Let $A$ be a subgroup of a group $G$. Then $A_{G}$ denotes the
maximal normal subgroup of $G$ contained in $A$. Let $G$ be a
group of order $p_1^{a_1}p_2^{a_2} \ldots p_k^{a_k}$, where
$p_1>p_2> \ldots >p_k$. We say that $G$ has an ordered Sylow tower
of supersolvable type if there exists a series
$$
1=G_0\leq G_1\leq G_2\leq \ldots \leq G_{k-1}\leq G_k=G
$$ of normal subgroups of $G$ such that for each $i=1,2,\ldots
,k$, $G_{i}/G_{i-1}$ is isomorphic to a Sylow $p_i$-subgroup of
$G$. By $G=[A]B$ we denote the semidirect product with normal
subgroup $A$ of $G$, $Z_n$ is a cyclic group of order $n$. We use
$d(G)$ to denote the derived length of a solvable group $G$.

Let $\mathfrak F$ and $\mathfrak H$ be non-empty formations. If
$G$ is a group then $G^{\mathfrak F}$ denotes the $\mathfrak
F$-residual of $G$, that is the intersection of all those normal
subgroups $N$ of $G$ for which $G/N \in \mathfrak F$. We define
$\mathfrak F \circ \mathfrak H =\{G \mid G^{\mathfrak H}\in
\mathfrak F\}$ and call $\mathfrak F \circ \mathfrak H$ the
formation product of $\mathfrak F$ and $\mathfrak H$ (see Doerk,
Hawkes, 1992, IV, 1.7). As usually, ${\mathfrak F}^2=\mathfrak F
\circ \mathfrak F$ and ${\mathfrak F}^n={\mathfrak F}^{n-1} \circ
\mathfrak F$ for every natural $n\geq 3$. A formation $\mathfrak
F$ is said to be saturated if $G/\Phi (G)\in \mathfrak F$ implies
that $G\in \mathfrak F$. In this paper, $\mathfrak N$ and
$\mathfrak A$ denotes the formations of all nilpotent and all
Abelian groups respectively. The other definitions and terminology
about formations could be referred to Doerk, Hawkes (1992),
Huppert (1967) and Shemetkov (1978).

{\bf    Lemma    2.1.}    {\sl Let $G$ be a bicyclic $p$-group.

$1.$ Let $N$ be a complemented normal subgroup in $G$. Then:

$1.1)$ if $p=2$, then $|N/\Phi(N)|\leq 4$;

$1.2)$ if $p>2$, then either $N=G$ or $N$ is cyclic.

$2.$ If $p>2$, then $G$ is metacyclic.

$3.$ If $p=2$, then any normal subgroup of $G$ is generated by at
most three elements.}

{\sl Proof.} 1. It follows from Monakhov, Gribovskaya (2001, Lemma
1).

2. It follows from Huppert (1967, III, 11.5).

3. Let $G=\langle a\rangle \langle b\rangle$ be a bicyclic
2-subgroup and $N$ a normal subgroup of $G$. Apply induction on
$|G|+|G/N|$. First we show that $|N/\Phi (N)|\le 8$. Assume that
$\Phi (N)\ne 1$. Then $\Phi (N)$ is normal in $G$ and by
induction, $N/\Phi (N)$ is generated by at most three elements.
Hence $|N/\Phi (N)|\le 8$ and by Huppert (1967, III, 3.15), $N$ is
generated by at most three elements. Consequently, $\Phi (N)=1$
and $N$ is an elementary Abelian group. By the inductive
assumption, $N$ is not contained in the proper bicyclic subgroups
of $G$. If $\langle a\rangle N\ne G$, then $\langle a\rangle N=
\langle a\rangle (\langle a\rangle N\cap \langle b\rangle)$ is
bicyclic, a contradiction. Hence $\langle a\rangle N=G$. Let
$T=\langle a\rangle \cap N$. Then $|T|\le 2$ and $G/T$ is bicyclic
2-subgroup with complemented normal subgroup $N/T$. By 1.1),
$|N/T|\le 4$. Hence $|N|\le 8$. The lemma is proved.

{\bf Example 2.2.} The calculations in the computer system GAP
(see GAP, 2009) show that the group $G$ of order $189=3^37$ having
number 7 in the library SmallGroups,
$$
G=<a,b,c,d \mid  b^3=c^3=d^7=1, \ a^3=c, \ [a,b]=c^{-1},
$$
$$
[a,d]=d^{-1}, \ [a,c]=[b,c]=[b,d]=[c,d]=1>,
$$
is the product of two cyclic subgroups $A=<bd>$ of order 21 and
$B=<ab>$ of order 9. Hence $G$ is bicyclic non-primary group of
odd order. There are only three non-trivial cyclic normal
subgroups in $G$: $N_1=<c>$ of order 3, $N_2=<d>$ of order 7,
$N_3=<cd>$ of order 21. Since $G/N_i$ is non-cyclic, it follows
that $G$ is non-metacyclic. Therefore the statement of Proposition
2 (Lemma 2.1) is not true for non-primary groups.

{\bf Example 2.3.}  The bicyclic 2-group $G$ of order 32,
$$ G=<a,b,c
\mid   a^2=b^8=c^2=1, \ [a,b]=c, \ [b,c]=b^4, \ [a,c]=1>,
$$
(see Huppert, 1953), contains a normal elementary Abelian subgroup
$N=<a>\times <b^4>\times {<c>}$ of order 8 with cyclic group $G/N$
of order 4. This example shows that the estimation of the number
of generators in Proposition 3 (Lemma 2.1) is exact.

Recall that $r_p(G)$ is the chief $p$-rank of the solvable group
$G$ (see Huppert, 1967, VI, 5.2). The chief rank is the maximum of
$r_p(G)$ for all $p\in \pi(G)$.

{\bf    Lemma    2.4.}    {\sl  Let $G$ be a solvable group having
a normal series such that every Sylow subgroup of its factors is
bicyclic. Then the orders of chief factors of $G$ are  $p$, $q^2$
or $8$, where $p$ and $q$ are primes from $\pi (G)$.}

{\sl Proof.} Let (1) be a normal series of $G$ such that every
Sylow subgroup of its factors is bicyclic. We refine this series
to a chief series of $G$. Let $\overline{N}=N/G_i$ be a minimal
normal subgroup of $\overline {G}=G/G_i$ such that $\overline{N}
\subseteq \overline{G_{i+1}}=G_{i+1}/G_i$.  Since $\overline {G}$
is solvable, $\overline{N}$ is an elementary Abelian $p$-subgroup
for some prime $p\in \pi (G)$. Besides, $\overline{N}$ is normal
in a bicyclic Sylow $p$-subgroup of $\overline{G_{i+1}}$. If
$p>2$, then $\overline{G_{i+1}}$ is metacyclic by Proposition~2
(Lemma 2.1). Hence $|\overline{N}|=p$ or $|\overline{N}|=p^2$. If
$p=2$, then $|\overline{N}|=2$, 4 or 8 by Proposition 3 (Lemma
2.1). As a result we obtain a chief series with factors of orders
$p$, $q^2$ or 8. By the Jordan-H\"older Theorem, all chief series
of some group are isomorphic. Hence $r_p(G)\leq 2$ for any prime
$p>2$ and $r_2(G)\leq 3$ by definition of the chief $p$-rank
$r_p(G)$. The lemma is proved.

{\bf    Lemma    2.5.}    {\sl Let $G$  be a group of odd order.
Then $G$ has a normal series such that every Sylow subgroup of its
factors is bicyclic if and only if the chief rank of $G$ is at
most 2.}

{\sl Proof.} Let $G$ has a normal series such that every Sylow
subgroup of its factors is bicyclic.  Then the chief rank of $G$
is at most 2 by Lemma 2.4. Conversely, if  the chief rank of $G$
is at most 2, then $G$ has a chief series in which every factor
either has prime order or is an elementary Abelian of order $p^2$
for some prime $p$. The lemma is proved.

{\bf    Lemma    2.6.}    {\sl Let $G$ be a solvable group having
a normal series such that every Sylow subgroup of its factors is
bicyclic. If $M$ is a maximal subgroup of $G$, then $|G:M|$ is
either a prime or the square of a prime or 8.}

{\sl Proof.} By Lemma 2.4, $G$ has a chief series
$$
1=G_0 < G_1 < \ldots < G_i < G_{i+1} < \ldots < G_m=G
$$
with factors of orders $p$, $q^2$ or $8$, where $p$ and $q$ are
primes. Let $G_i\subseteq M$, but $G_{i+1}\not \subseteq M$. Since
$M$ is maximal in $G$, it follows that $G_{i+1}M=G$ and $
|G:M|=|G_{i+1}:G_{i+1}\cap M|. $ Because $G_i\subseteq G_{i+1}\cap
M$, we have
$$
|G_{i+1}:G_{i+1}\cap M|=\frac{|G_{i+1}:G_i|}{|G_{i+1}\cap
M:G_i|}$$ and $|G:M|$ is either a prime or the square of a prime
or 8. The lemma is proved.

\medskip

{\bf    Lemma    2.7.}  (Bloom, 1967, Theorem 3.4) {\sl Let $G$ be
a subgroup of $GL(2,q)$ and $q=p^\alpha$, where $p$ is prime.
Then, up to conjugacy in $GL(2,q)$, one of the following occurs:

$1)$ $G$ is cyclic;

$2)$ $G=QM$, where $Q$ is a subgroup of the $p$-group $\left\{\left(%
\begin{array}{cc}
  1 & 0 \\
  \tau & 1 \\
\end{array}%
\right)\mid \tau \in GF(q)\right\}$ and $M\subseteq N_G(Q)$ is a
subgroup of the group $D$ of all diagonal matrices;

$3)$ $G=\{Z_u,S\}$, where $u$ divides $q^2-1$, $S:Y\rightarrow
Y^q$, for all $Y\in Z_u$, and $S^2$ is a scalar 2-element in
$Z_u$;

$4)$ $G=\{M,S\}$, where $M\subseteq D$ and  $|G:M|=2$;

$5)$ $G=\langle SL(2,p^\beta),V\rangle$ ("Case 1") or
$$
G=\Big{\langle} SL(2,p^\beta),V, \left(
\begin{array}{cc}
  b & 0 \\
  0 & \epsilon b \\
\end{array}
\right)\Big{\rangle} ,
$$
("Case 2"), where $V$ is a scalar matrix, $\epsilon$ generates
$(GF(p^\beta))^{*}$, $p^\beta>3$, $\beta | \alpha$. In Case 2,
$|G:\langle SL(2,p^\beta),V\rangle|=2$;

$6)$ $G/\{-I\}$ is isomorphic to $S_4\times Z_u$, $A_4\times Z_u$
or $A_5\times Z_u$, if $p\neq 5$, where $Z_u$ is a scalar subgroup
of $GL(2,q)/\{-I\}$;

$7)$ $G$ is not of type $(6)$, but $G/\{-I\}$ contains $A_4\times
Z_u$ as a subgroup of index 2, and $A_4$ as a subgroup with cyclic
quotient group, $Z_u$ is as in type $(6)$ with $u$ even.}

{\bf    Lemma    2.8.}    {\sl Let $H$ be an $A_4$-free
$p^{\prime}$-subgroup of $GL(2,p)$, where $p$ is prime. Then $H$
is metabelian.}

{\sl Proof.} We shall use the result of Lemma 2.7. A subgroup $H$
from Proposition 1 is Abelian. The order of a subgroup $H$ from
Proposition 2 is divisible by a prime $p$. Since the group of all
diagonal matrices is Abelian, it follows that a subgroup $H$ from
Proposition 3-4 is metabelian. A subgroup $H$ from Proposition 5-7
is not $A_4$-free. Hence if $H$ is an $A_4$-free $p'$-subgroup
$GL(2,p)$, then $H$ is metabelian. The lemma is proved.

{\bf    Lemma    2.9.}    {\sl Let $H$ be a subgroup of $GL(3,2)$.
Then $H\in \{1$, $GL(3,2)$,  $Z_{2}$,  $Z_{3}$, $Z_{7}$,
 $Z_{2}\times Z_{2}$, $Z_{4}$, $D_{8}$, $S_{3}$, $A_{4}$, $S_{4}$, $[Z_{7}]Z_{3}\}$.}

{\sl Proof.} By  Huppert (1967, II, 6.14), $GL(3,2)\simeq
PSL(2,7)$. In view of Huppert (1967, II, 8.27), we conclude that
$H$ satisfies the hypotheses of our lemma.

{\bf    Lemma    2.10.}    {\sl Let $G$ be a solvable group such
that the index of each of its maximal subgroup is either a prime
or the square of a prime or 8. Then the following statements hold:

$1)$ $G\in\mathfrak{N}_{2'}\circ \mathfrak{N}_{2}\circ
\mathfrak{U}$. In particular, the nilpotent length of $G$ is at
most 4;

$2)$ $G$ contains a normal subgroup $N$ such that $G/N$ is
supersolvable and $N$ possesses an ordered Sylow tower of
supersolvable type;

$3)$ $l_{2} (G)\le 2$, $l_{3} (G)\le 2$ and $l_{p}(G) \le 1$ for
every prime $p>3$. If $G$ is a group of odd order, then $l_{p}(G)
\le 1$ for every prime $p\in \pi(G)$;

$4)$ $G$ contains a normal Hall
$\{2,3,7\}^{\prime}$\nobreakdash-\hspace{0pt}subgroup $H$ and $H$
possesses an ordered Sylow tower of supersolvable type;

$5)$ if $G$ is a group of odd order, then $G$ possesses an ordered
Sylow tower of supersolvable type.}

{\sl Proof.} 1. It follows from Gribovskaya (2001, Theorem 2,
Corollary 3).

2. By 1) $G\in\mathfrak{N}_{2'}\circ \mathfrak{N}_{2}\circ
\mathfrak{U}$, i.e. $G^{\mathfrak{U}}\in\mathfrak{N}_{2'}\circ
\mathfrak{N}_{2}$. Hence $G^{\mathfrak{U}}=[T]H$, where $T$ is a
$2'$--Hall subgroup, $H$ is a Sylow $2$-subgroup. Since $T\in
\mathfrak{N}_{2'}$, it follows that $T$ is nilpotent and
$G^{\mathfrak{U}}$ possesses an ordered Sylow tower of
supersolvable type.

3. We use induction on $|G|$. Let $p$ be a prime divisor of $|G|$.
By Huppert (1967, VI, 6.9), we may assume that $O_{p'} (G)=\Phi
(G)=1$ and $G=[F]M$, where the Fitting subgroup $F=F(G)=C_G(F)$ is
the unique minimal normal $p$-subgroup and $M$ is a maximal
subgroup of $G$. Hence a Sylow $p$-subgroup $G_{p} =[F](G_{p} \cap
M)=[F]M_{p},$ where $M_{p} $ is a Sylow $p$-subgroup of $M$. If
$M_{p} =1$, then $F=G_{p} $ and $l_{p} (G)\le 1$. Let $M_{p} \ne
1$. Since $|F|=|G:M|$, it follows that $|F|$ is equal either to
$p$ or $p^{2}$, or 8. If $|F|=p$, then $G/F$ is a cyclic group
whose order divides $(p-1)$. Hence $G_{p} =F$, a contradiction.

Let $|F|=p^{2} $. Then $G/F$ is isomorphic to a subgroup of
$GL(2,p)$. Since $|GL(2,p)|=(p^{2} -p)(p^{2} -1)$, the order of
$G_{p} $ is equal to $p^{3}$ and by Huppert (1967, VI, 6.6),
$l_{p} (G)\le 2$. Since $F=C_G(F)$, $G_{p} $ is non-Abelian and by
Huppert (1967, I, 14.10), it is isomorphic either to a metacyclic
group $M_{3}(p)=\langle a,b \mid a^{p^{2}}=b^{p}=1, a^{b}=a^{1+p}
\rangle =[\langle a\rangle ]\langle b\rangle ,$ or to a group of
exponent $p$. Since $\Omega _{1} (M_{3} (p))$ is an elementary
Abelian $p$-subgroup of order $p^{2} $, it hasn't the complement
in $M_{3} (p)$. Hence $G_{p}$ is a group of exponent $p$. If $G$
has a odd order or $p$ is not a Fermat prime, then by Huppert,
Blackburn~(1982, IX, 4.8), $l_{p} (G)\le 1$. But now by Huppert,
Blackburn~(1982, IX, 5.5(b)), $l_{p} (G)\le 1$ for $p>3$.

Finally, let $|F|=8$. Then $p=2$ and $G/F$ is isomorphic to a
subgroup $H$ of $GL(3,2)$. In this case $O_{2} (G/F)=1$ and by
Lemma 2.9, $H\in \{ Z_{3} ,\, Z_{7} ,\, S_{3} ,\, [Z_{7} ]Z_{3}
\}$. Evidently, $l_{2} (G)\le 2$.

4. We show that $G$ has a normal Hall $\pi$--subgroup $G_{\pi}$
for $\pi=\pi(G) \setminus \{2,3,7\}$. Since the class of all
$\pi$-closed subgroups is a saturated formation, by induction we
can assume that $O_{\pi}(G)=1$ and the Fitting subgroup $F$ is an
elementary Abelian $p$-subgroup whose order divides $2^{3}$,
$3^{2}$ or $7^{2}$. Hence the group $G/F$ is isomorphic to a
subgroup of $GL(n,p)$ for $p=2$ and $n \le 3$, or for $p \in
\{3,7\}$ and $n \le 2.$ Since $\pi(GL(n,p)) \subseteq \{2,3,7\}$
for given $n$ and $p$, it follows that $G$ is a $\pi'$--subgroup.

By Monakhov, Selkin, Gribovskaya (2002, Corollary 2.4), $G_{\pi}$
possesses an ordered Sylow tower of supersolvable type.

5. It follows from Monakhov, Selkin, Gribovskaya (2002, Corollary
2.3).

\section{PROOFS OF THEOREM 1.1 AND COROLLARY 1.2-1.3}

{\bf Proof of Theorem 1.1}

By Lemma 2.6 and Lemma 2.10 (1-4), we must only prove that the
derived length of $G/\Phi(G)$ is at most 5.

We first show that $G\in \mathfrak{N}\circ \mathfrak{A}^4$. Apply
induction on $|G|$. Assume that $\Phi(G)\not = 1$.  Since any
quotient group satisfies the hypothesis of the theorem,
$G/\Phi(G)\in \mathfrak{N}\circ \mathfrak{A}^4$ by induction.
Since $\mathfrak{N}\circ \mathfrak{A}^4$ is a saturated formation,
it follows that $G\in \mathfrak{N}\circ \mathfrak{A}^4$. Next we
assume that $\Phi(G)=1$.

Now suppose that the Fitting subgroup $F(G)$ is not a minimal
normal subgroup in $G$. Then $F(G)$ is the direct product of
minimal normal subgroups of $G$, i.\,e. $F(G)=F_{1}\times F_{2}
\times \ldots \times F_{n}$, where $F_{i}$ is a minimal normal
subgroup of $G$ for any $i$ and $n \geq 2$. By the inductive
assumption, we have $G/F_{i}\in \mathfrak{N}\circ \mathfrak{A}^4$.
Consequently, $G\in \mathfrak{N}\circ \mathfrak{A}^4$, because
$\mathfrak{N}\circ \mathfrak{A}^4$ is a formation.

Next we assume that $F=F(G)$ is the unique minimal normal subgroup
of $G$. Besides, $F = C_{G}(F)$ and $G = [F]M$, where $M$ is a
maximal subgroup of $G$. Since $|F|=|G:M|$, it follows by Lemma
2.6 that $|F|$ is equal to $p$, $p^{2}$ or 8, where $p$ is prime.

If $|F| = p$, then $G/F$ is cyclic, since it is the subgroup of
$\mathrm{Aut} F=Z_{p-1}$. Hence $G/F \in \mathfrak{A}$. Let $|F| =
p^{2}$. Then $G/F$ is isomorphic to an irreducible solvable
subgroup of $GL(2,p)$. By Monakhov, Gribovskaya (2001, Lemma 3),
$G/F\in \mathfrak{A}^4$.

It remains to study the case $|F| = 8$. Then $G/F$ is isomorphic
to a solvable subgroup $H$  of $GL(3,2)$. Let's notice that $F$ is
the maximal normal 2-subgroup of $G$, i.\,e $F = O_{2} (G)$. Hence
$O_{2} (G/F) = 1$. By Lemma 2.9,  $G/F \in
\{Z_{3},S_3,Z_{7},[Z_{7}]Z_{3}\}$ and $G/F \in \mathfrak{A}^2
\subseteq \mathfrak{A}^4$.

From all the above, we proved that $G/F\in \mathfrak{A}^4$. As $F$
is nilpotent, $G \in \mathfrak{N}\circ \mathfrak{A}^4$. Since
$F/\Phi(G)$ is Abelian and $(G/\Phi(G))/(F/\Phi(G))\simeq G/F$, it
follows that $G/\Phi(G)\in \mathfrak{A}^5$ and $d(G/\Phi(G))\leq
5$. The theorem is proved.

\medskip

{\bf Proof of Corollary 1.2}

1. By Proposition 3 (Theorem 1.1), we obtain $l_2(G)\leq 2$,
$l_3(G)\leq 2$ and $l_{p}(G) \le 1$ for every prime $p>3$. Now we
show that $l_p(G)\leq 1$, where $p\in \{ 2,3\}$. By Huppert (1967,
VI, 6.9), we may say that $O_{p'}(G)=\Phi(G)= 1$. By Lemma 2.4,
the Fitting subgroup $F=F(G)$ is the unique minimal normal
subgroup of order $p^\alpha$, where $\alpha\leq 3$ for $p=2$ and
$\alpha\leq 2$ for $p=3$. In particular, $C_G(F)=F$ and $G=[F]M$
for some maximal subgroup $M$ of $G$. If $|F|=p$, then $G/F$ is
isomorphic to a subgroup of order $p-1$ and $l_{p}(G)\le 1$. If
$|F|=4$, then $\mbox{Aut}(F(G))\simeq GL(2,2)\simeq S_3$. Hence
either $G/F(G)\simeq Z_3$ or $G/F(G)\simeq S_3$. If $G/F(G)\simeq
Z_3$, then $G\simeq A_4$. If $G/F(G)\simeq S_3$, then $G\simeq
S_4$. It means that $G$ is not $A_4$-free, a contradiction.

Now let $|F|=8$. Then $G/F$ is isomorphic to a subgroup of
$GL(3,2)$. Since $O_2(G/F)=1$, it follows by Lemma 2.9, that
$G/F\in \{Z_{3},S_3,Z_{7},[Z_{7}]Z_{3}\}$. In all cases, except
$G/F\simeq S_3$, we have $l_{2}(G)\le 1$. Suppose that $G/F$ is
isomorphic to $S_3$. We may construct the subgroup $H=[F]Z_3$ in
$G$. Then the alternating group $A_4$ of degree 4 is contained in
$H$, a contradiction.

Let $|F|=9$. Then $G/F$ is isomorphic to a subgroup of $GL(2,3)$
and $O_3(G/F)=1$. It is well known that $H\in \{ 1$, $Z_{2} $,
$Z_{4} $, $Z_{8} $, $Z_{2} \times Z_{2} $, $D_{8} $, $Q_{8} $,
$SD_{16} $, $SL(2,3)$, $GL(2,3)\} $. In any case, except $G/F\cong
SL(2,3)$ and $G/F\cong GL(2,3)$, $F$ is a Sylow 3-subgroup in $G$
and $l_{3} (G)\le 1$. Since $SL(2,3)$ and $GL(2,3)$ are not
$A_{4}$-free, we have a contradiction.

2. We use induction on $|G|$. We first prove that $G\in
\mathfrak{N}\circ \mathfrak{A}^2$. By induction, we can assume
that $\Phi (G)=1$ and $G$ has the unique minimal normal subgroup
which coincides with Fitting subgroup $F=F(G)$. By Proposition 1
(Corollary 1.2), $l_p(G)\leq 1$. Hence $F$ is a Sylow $p$-subgroup
of $G$. Besides, $F = C_{G}(F)$ and $F$ has a complement $M$ in
$G$, where $M$ is a maximal subgroup of $G$. By Lemma 2.6, $|F|$
is equal to $p$, $p^{2}$ or 8, where $p$ is prime.

If $|F| = p$, then $G/F$  is cyclic, since it is the subgroup of
$\mathrm{Aut} F=Z_{p-1}$. Hence $G/F$ is Abelian. Let $|F| =
p^{2}$. Then $G/F$ is isomorphic to an irreducible solvable
$p^{\prime}$-subgroup $H$ of $GL(2,p)$. By Lemma 2.8, $H$ is
metabelian, i.\,e. $G/F\in \mathfrak{A}^2$.

Now let $|F|=8$. Then $G/F$ is isomorphic to a subgroup of
$GL(3,2)$. By Lemma 2.9, $G/F\in \{Z_{3},Z_{7},[Z_{7}]Z_{3}\}$.
Then $H$ is metabelian and $G/F\in \mathfrak{A}^2$.

So, in any case $G/F\in \mathfrak{A}^2$. Since $F/\Phi(G)$ is
Abelian and $(G/\Phi(G))/(F/\Phi(G))\simeq G/F$, it follows that
$G/\Phi(G)\in \mathfrak{A}^3$ and $d(G/\Phi(G))\leq 3$. The
corollary is proved.

\medskip

{\bf Proof of Corollary 1.3}

1. By Lemma 2.10 (5), our assertion holds.

2. We show that the {derived subgroup} of $G$ is nilpotent. We use
induction on $|G|$. Without loss of generality, we may assume that
$\Phi(G)=1$ and  $G$ has a unique minimal normal subgroup which
coincides with Fitting subgroup $F=F(G)$. Then $F$ is an
elementary Abelian $p$-subgroup for some prime $p$. Since
$\Phi(G)=1$, it follows that $G$ has a maximal subgroup $M$ such
that $G = [F]M$. Because $|F|=|G:M|$, we have by Lemma~2.6, that
$|F|$ is equal to $p$ or $p^{2}$. By Proposition 3 (Lemma 2.10),
$l_p(G)=1$. Hence $F$ is a Sylow $p$-subgroup of $G$ and $G/F$ is
a $p^{\prime}$-subgroup. In the solvable groups the Fitting
subgroup coincides with its centralizer in $G$, hence $G/F$ is
isomorphic to a subgroup of $\mathrm{Aut} F$.

If $|F| = p$, then $G/F$ is cyclic and $G^{\prime}\subseteq F$.
Let $|F| = p^{2}$. Then $G/F$ is isomorphic to an irreducible
solvable $p^{\prime}$-subgroup $H$ of $GL(2,p)$. By Dixon (1971,
Theorem 5.2), $H$ is Abelian and $G^{\prime}\subseteq F$. So, in
any case the derived subgroup of $G$ is nilpotent.

Since $F/\Phi (G)$ is Abelian, it follows that $G/\Phi (G)$ is
metabelian. The corollary is proved.

\medskip

{\bf Example 3.1.} Let $E_{7^2}$ be an elementary Abelian group of
order $7^2$. The automorphism group of $E_{7^2}$ is the general
linear group $GL(2,7)$ with cyclic center $Z=Z(GL(2,7))$ of order
$6$. We choose a subgroup $C$ of order $2$ in $Z$. Evidently, $C$
is normal in $GL(2,7)$. The calculations in the computer system
GAP show that $GL(2,7)$ has a subgroup $S$  of order $48$ such
that $S/C$ is isomorphic to the symmetric group $S_4$ of degree
$4$. The semidirect product $G=[E_{7^2}]S$ is a group of order
$2352=2^47^23$. In particular, $\Phi(G)=1$. The nilpotent length
of $G$ is equal to $4$, the derived length of $G$ is equal to $5$.
The group $G$ has the chief series
$$
1\subset E_{7^2}\subset [E_{7^2}]Z_2\subset [E_{7^2}]Q_8\subset
[[E_{7^2}]Q_8]Z_3\subset [E_{7^2}]S=G
$$
with bicyclic factors:
$$
E_{7^2}, \  ([E_{7^2}]Z_2)/(E_{7^2})\simeq Z_2, \
([E_{7^2}]Q_8)/([E_{7^2}]Z_2)\simeq E_4,
$$
$$
([[E_{7^2}]Q_8]Z_3)/([E_{7^2}]Q_8)\simeq Z_3, \
(G/[[E_{7^2}]Q_8]Z_3)\simeq Z_2.
$$
Hence the estimations of the nilpotent length and the derived
length, which are obtained in Theorem 1.1, are exact.

{\bf Example 3.2.} Let $E_{5^2}$ be an elementary Abelian group of
order $5^2$. The automorphism group of $E_{5^2}$ is the general
linear group $GL(2,5)$. The group $GL(2,5)$ has a subgroup, which
is isomorphic to the symmetric group $S_3$ of degree 3. The
semidirect product $G=[E_{5^2}]S_3$ is an $A_4$-free group with
identity Frattini subgroup. The derived length of $G$ is equal to
3. The group $G$ has the chief series
$$
1\subset E_{5^2}\subset [E_{5^2}]Z_3\subset [E_{5^2}]S_3=G
$$
with bicyclic factors:
$$
E_{5^2}, \  ([E_{5^2}]Z_3)/(E_{5^2})\simeq Z_3, \
([E_{5^2}]S_3)/([E_{5^2}]Z_3)\simeq Z_2.
$$
Consequently, the estimation of the derived length, which is
obtained in Corollary 1.2, is exact.

{\bf Example 3.3.} It is well known that $S_4$ has the normal
series
$$1\leq E_4\leq A_4\leq S_4$$ with bicyclic factors and
$l_{2}(S_4)=2$. The group $G=[E_{3^2}]SL(2,3)$ has the normal
series
$$1\leq E_{3^2}\leq [E_{3^2}]Z_2\leq [E_{3^2}]Q_8\leq
[E_{3^2}]SL(2,3)$$ with bicyclic factors and $l_{3}(G)=2$.

The project is supported by the Belarus republican fund of basic
researches (No. F 08R-230 ).

\section{REFERENCES}

\ \ \ \ Bloom, D. (1967) The subgroups of $PSL(3,q)$ for odd $q$.
{\sl Trans. Amer. Math. Soc.} 1(127):150--178.

Dixon, J. D. (1971) {\sl The structure of linear groups.} Van
Nostrand, Princeton, N. J. and London.

Doerk, K., Hawkes, T. (1992) {\sl Finite soluble groups.} Berlin,
New York: Walter de Gruyter.

GAP (2009) Groups, Algorithms, and Programming, Version 4.4.12.
{\sl www.gap-system.org.}

Gribovskaya, E. E. (2001) Finite solvable groups with the index of
maximal subgroups is $p$, $p^{2}$ or 8.  {\sl Vesti NAN Belarus.}
4:11--14. [In Russian]

Huppert, B. (1953) Uber das Produkt von paarweise vertauschbaren
zyklischen Gruppen. {\sl Math. Z.} 58:243--264.

Huppert, B. (1967) {\sl Endliche Gruppen I.} Berlin, Heidelberg,
New York: Springer.

Huppert, B., Blackburn, N. (1982) {\sl Finite Groups II.} Berlin,
Heidelberg, New York: Springer.

Monakhov, V. S., Gribovskaya, E. E. (2001) Maximal and Sylow
subgroups of solvable finite groups. {\sl Matem. Notes}.
70(4):545--552.

Monakhov, V. S., Selkin , M. V., Gribovskaya, E. E. (2002) On
normal solvable subgroups of finite groups. {\sl Ukr. Math. J.}
54(7):950--960. [in Russian]

Shemetkov, L. A. (1978) {\sl Formations of finite groups.}  M.:
Nauka. [in Russian]

\bigskip

Department of Mathematics, Gomel Francisk Skorina State
University, Gomel 246019, Belarus; e-mail: monakhov@gsu.by,
e-mail: trofim08@yandex.ru

\end{document}